\documentclass[12pt]{amsart}
\usepackage[utf8]{inputenc}
\usepackage{amsmath,amsthm}
\usepackage{comment}
\usepackage{amssymb}
\usepackage{stmaryrd}
\usepackage{hyperref}
\usepackage{color}
\usepackage[a4paper,left=3cm,right=3cm,top=3cm,bottom=4.5cm]{geometry}
\usepackage{pdfsync}
\usepackage{tikz-cd}
\usepackage{url}
\usepackage{xpatch}
\xpatchcmd{\paragraph}{\normalfont}{{\normalfont\bfseries}}{}{}

\newcommand{\defterm}[1]{\emph{#1}}

\theoremstyle{plain}
	\newtheorem{thm}{Theorem}
	\numberwithin{thm}{section}
	\newtheorem*{thm*}{Theorem}
	\newtheorem{cor}[thm]{Corollary}
	\newtheorem*{cor*}{Corollary}
	
	\newtheorem*{prop*}{Proposition}
	\newtheorem{lem}[thm]{Lemma}
	\newtheorem*{lem*}{Lemma}
	
	\newtheorem*{ex*}{Exercise}
	
	\newtheorem*{claim*}{Claim}
	
	\newtheorem*{conj*}{Conjecture}
	\newtheorem{question}[thm]{Question}
	\newtheorem*{question*}{Question}
	
	\newtheorem*{obs*}{Observation}

\theoremstyle{definition}
	\newtheorem{Def}[thm]{Definition}
	\newtheorem*{Def*}{Definition}
	\newtheorem{rmk}[thm]{Remark}
	\newtheorem*{rmk*}{Remark}
	
	\newtheorem{soln*}{Solution}
	
	\newtheorem*{note*}{Note}
	\newtheorem{eg}[thm]{Example}
	\newtheorem*{eg*}{Example}	
	
	\newtheorem*{construction*}{Construction}
	
	\newtheorem*{warning*}{Warning}
	\newtheorem{fact}[thm]{Fact}
	\newtheorem*{fact*}{Fact}	
	
\newcommand{\ints}{\mathbb{Z}}

\newcommand{\reals}{\mathbb{R}}

\newcommand{\nats}{\mathbb{N}}
\newcommand{\pjv}{\mathbb{P}}

\newcommand{\id}{\mathrm{id}}

\newcommand{\Hom}{\mathrm{Hom}}
\DeclareMathOperator{\dom}{\mathrm{dom}}
\DeclareMathOperator{\cod}{\mathrm{cod}}
\newcommand{\Ob}{\operatorname{Ob}}
\newcommand{\Mor}{\operatorname{Mor}}

\newcommand{\Psh}{\mathsf{Psh}}

\newcommand{\Cat}{\mathsf{Cat}}

\newcommand{\bbB}{\mathbb{B}}

\newcommand{\calA}{\mathcal{A}}

\newcommand{\calC}{\mathcal{C}}
\newcommand{\calD}{\mathcal{D}}

\newcommand{\nerve}{\mathbf{N}}

\newcommand{\Gal}{\mathrm{Gal}}

\newcommand{\Modhom}{\mathsf{Mod}^\mathrm{hom}}
\newcommand{\Modstr}{\mathsf{Mod}^\mathrm{str}}
\newcommand{\Modelem}[1]{\mathsf{Mod}^\mathrm{elem_{#1}}}
\newcommand{\Str}{\mathsf{Str}}

\newcommand{\mono}{\mathrm{mono}}
\newcommand{\Ext}{\mathrm{Ext}}
\newcommand{\Lift}{\mathrm{Lift}}
\newcommand{\sd}{\mathrm{sd}}
\newcommand{\Asph}{\mathrm{Asph}}
\newcommand{\Inj}{\mathrm{Inj}}
\newcommand{\contr}{\mathrm{contr}}

\newcommand{\sSet}{\mathsf{sSet}}

\newcommand{\bigdoublevee}{%
  \mathop{
    \mathchoice{\bigvee\mkern-15mu\bigvee}
               {\bigvee\mkern-12.5mu\bigvee}
               {\bigvee\mkern-12.5mu\bigvee}
               {\bigvee\mkern-11mu\bigvee}
    }
}

\newcommand{\true}{\top}
\newcommand{\false}{\bot}

\title{Homotopy Types of Abstract Elementary Classes}
\author[Tim Campion]{Tim Campion}
\address{University of Notre Dame\\
Department of Mathematics\\
255 Hurley, Notre Dame\\
IN 46556, USA.
}
\email{tcampion@nd.edu}
\author[Jinhe Ye]{Jinhe Ye}
\address{University of Notre Dame\\
Department of Mathematics\\
255 Hurley, Notre Dame\\
IN 46556, USA.
}
\email{jye@nd.edu}
\date{September 2019}
\begin{document}
\maketitle
\begin{abstract}
   We prove that for any homotopy type $X$, there is an abstract elementary class $\calC$, with joint embedding, almagamation and no maximal models such that the classifying space realizes the homotopy type $X$. We provide a few explicit examples.
  \end{abstract}
\section{Introduction}
In \cite{classify}, given a first order theory $T$, we studied the homotopy type of $|\Modelem\omega(T)|$, the classifying space of the category of models of $T$ with elementary embeddings as morphisms. In particular, in all the examples where we can determine its homotopy type, the higher homotopy always vanishes. It leads us to the following question.
\begin{question}\label{question:hom}
\begin{enumerate} Let $X$ a homotopy type given.
    \item Is there a first order theory $T$ such that $|\Modelem\omega(T)|$  realizes the homotopy type $X$?
    \item  Is there an abstract elementary class $\calA$ with joint embedding, almagamation and no maximal models such that the classifying space of $\calA$ realizes the homotopy type $X$?
\end{enumerate} 
\end{question}
In \cite{classify}, a partial answer to (2) was given. In this paper, we answer Question \ref{question:hom}(2) positively. In particular, for any small category $\calC$ of monomorphisms, we consider the category $\calA$ of presheaves $X$ on $\calC$ such that the slice category $\calC \downarrow X$ is weakly contractible, and satisfying the technical condition that for each representable presheaf $C$, each map $C \to X$ is a monomorphism. We show (Theorem \ref{thm:mainthm}) that $\calA$ is an AEC with amalgamation which typically has no maximal models, and that $\calA$ has the same homotopy type as $\calC$. Because every homotopy type is realized by a poset, and in particular a category of monomorphisms \cite{thomason}, the conclusion follows. In fact, $\calA$ is axiomatized via a basic theory in $L_{\kappa,\omega}$, where $\kappa = \max(\aleph_0, \#\calC)$ (here $\#\calC$ is the number of morphisms of $\calC$). The proof uses a well-known folklore characterization of weakly contractible simplicial sets via a lifting condition (Fact \ref{fact:contr}).

We conclude in Section \ref{sec:examples} by illustrating this example and a close variant, describing explicit examples of AECs with certain interesting homotopy types. including $S^1$ (Example \ref{eg:dirtree}) $\reals \pjv^\infty$ (Example \ref{eg:undirtree}), and $S^2$ (Example \ref{eg:s2}). We believe that such examples may be of interest in their own right, quite apart from considerations from a mathematical logic perspective.
\paragraph{Conventions}
We fix a strongly inaccessible cardinal $\lambda$. For each ordinal $\alpha$, let $V_\alpha$ denote the sets of rank less than $\alpha$. Sets in $V_{\lambda}$ are called \defterm{small}. The large cardinal hypothesis only enters the picture for convenience, our results will continue to hold, \emph{mutatis mutandis}, in ZFC.\\
Unlike in \cite{classify}, we work with $\nerve\calC$ and simplicial homotopy theory instead of geometric realizations and homotopy theory of topological spaces. For readers unfamiliar with simplicial homotopy theory, we recommend \cite{goerss_jardine}.

\section{Backgrounds in simplicial homotopy theory}
Let $\Delta$ denote the simplex category, i.e. the category of finite nonempty ordinals $[n] = \{0,1,\dots, n\}$ and order-preserving maps. Recall that the category $\sSet$ of \defterm{simplicial sets} is the category of presheaves on $\Delta$. We denote the representable on $[n]$ by $\Delta[n]$. We denote by $\partial \Delta[n] \subseteq \Delta[n]$ the subpresheaf obtained by deleting the top cell, so that its set of $m$-simplices $(\partial \Delta[n])_m$ is the set of non-epimorphic morphisms $[m] \to [n]$. We write $i: \partial \Delta[n] \to \Delta[n]$ for the canonical inclusion map.

In the following, we will often silently identify a category with its nerve.

\begin{Def}[Barycentric Subdivision]
For $[n] \in \Delta$, let $\sd[n]$ denote the nerve of the poset of nonempty subsets of $[n]$. Let $p: \sd[n] \to \Delta [n]$ denote the map which sends a set to its maximum element. This construction is functorial, and by left Kan extension we obtain a functor $\sd: \sSet \to \sSet$, the \defterm{barycentric subdivision functor}, and a natural transformation $\sd \Rightarrow \id$. We denote the composite $p^l = p (\sd p) \dots (\sd^{l-1} p)$.
\end{Def}

Note that if $X$ is a simplicial complex, then $\sd X$ is a simplicial complex and (the nerve of) a poset. In particular, the iterated barycentric subdivisions $\sd^k \Delta[n]$ (for $k \geq 0$) and $\sd^k \partial \Delta[n]$ (for $k \geq 1$) are posets.

\begin{fact}[\cite{bourke}, Section 3.2]\label{fact:contr}
A simplicial set $X$ is weakly contractible if and only if the following condition is satisfied. For every $n \in \nats$ and $k \geq 1$, and every map $f: \sd^k(\partial \Delta[n]) \to X$, there exists $l \in \nats$ and a map $g: \sd^{k+l}(\Delta[n]) \to X$ such that $g(\sd^{k+l} i) = f p^l$. Here $i: \partial \Delta[n] \to \Delta[n]$ is the canonical inclusion as above, and $p^l: \sd^{k+l} (\partial \Delta[n]) \to \sd^k (\partial \Delta[n])$ is the collapse map as above.

In particular, a category $\calC$ is weakly contractible if and only if the following condition is satisfied. For every $n \in \nats$ and $k \geq 1$, and every functor $F: \sd^k(\partial \Delta[n]) \to \calC$, there exists $l \in \nats$ and a functor $G: \sd^{k+l}(\Delta[n]) \to X$ such that $G(\sd^{k+l} i) = F p^l$.
\end{fact}

\section{Abstract Elementary Classes and Infinitary Logic}

We work with infinitary first-order logic. We always work over a (small) finitary signature, i.e. all atomic function and relation symbols have finite arities. We allow our signatures to have arbitrarily many sorts. We write $L_{\kappa,\lambda}$ for infinitary first-order logic where conjunctions and disjunctions are limited to be over sets of formulae of cardinality $<\kappa$ and quantifiers range over tuples of variables of length $<\lambda$, where $\kappa,\lambda$ are either infinite cardinals or else $\infty$, indicating that there is no restriction on the size of conjunction / disjunction or quantification, as appropriate. We write $L_\kappa$ for $L_{\kappa,\kappa}$. In particular, full infinitary first-order logic is denoted by $L_\infty$. For a signature $\Sigma$, we use $L_{\kappa,\lambda}(\Sigma)$ denote the $L_{\kappa,\lambda}$-formulae in the signature $\Sigma$.  Note that among the sentences of $L_{\kappa,\lambda}(\Sigma)$ are the empty conjuction $\true$ and the empty disjunction $\false$. A \defterm{theory} $T$ in  $L_{\kappa,\lambda}(\Sigma)$ is a (small) set of sentences of $L_{\kappa,\lambda}(\Sigma)$.

\begin{Def}[Basic theory {\cite[Definition 5.31]{lpacc}}]
Let $\Sigma$ be a finitary signature over some set of sorts. A formula in $L_\infty(\Sigma)$ is called
\begin{enumerate}
    \item \defterm{positive-primitive} if it is of the form $\exists y (\psi(x,y))$ where $\psi(x,y)$ is a (possibly infinite) conjuction of atomic formulae (here $x,y$ are possibly infinite tuples);
    \item \defterm{positive-existential} if it is a (possibly infinitary) disjunction of positive-primitive formulae;
    \item \defterm{basic} if is of the form $\forall x (\varphi(x) \to \psi(x))$ where $\varphi(x),\psi(x)$ are positive-existential.
\end{enumerate}
A theory $T$ of $L_\infty$ is called \defterm{basic} if it consists of basic sentences.
\end{Def}

\begin{eg}
Any \defterm{universal-positive} formula, i.e. a formula of the form $\forall x ( \phi(x,y))$ where $\phi(x,y)$ is a (possibly infinite) conjunction of atomic formulae, is basic. Any \defterm{universal-positive} theory, i.e. a theory consisting only universal-positive sentences, is basic.
\end{eg}
\begin{Def}\label{def:hom-str-elem}
Let $\Sigma$ be a finitary signature over some set of sorts and let $T$ be a theory in $L_\infty(\Sigma)$. Let $M,N$ be two models of $T$, a map $\phi:M\to N$ is
\begin{enumerate}
    \item a \defterm{homomorphism} if for any relation $R\in \Sigma$, function $f\in \Sigma$ and $m \in M$ a tuple of appropriate length, $R(m)$ holds in $M$ implies $R(\phi(m))$ holds in $N$ and for any function $f\in \Sigma$, $f(\phi(m))=\phi(f(m))$.
    \item a \defterm{strong embedding} if $\phi$ is a embedding and homomorphism such that for any relation $R\in \Sigma \cup \{=\}$ and $m\in M$ a tuple of appropriate length, $\neg R(m)$ holds in $M$ implies $\neg R(\phi(m))$ holds in $N$.
    \item a \defterm{$L_{\kappa,\lambda}$-elementary embedding}  if $\phi$ is a strong embedding such that for any $\psi(x)\in L_{\kappa,\lambda}(\Sigma)$ and a tuple $m\in M$ of appropriate length, $\psi(m)$ holds in $M$ iff $\psi(\phi(m))$ holds in $N$.
\end{enumerate}

\end{Def}
\begin{Def}[$\Modhom(T),\Modstr(T),\Modelem\lambda(T)$]\label{def:mod}
Let $\Sigma$ be a finitary signature over some set of sorts and let $T$ be a theory of $L_\infty(\Sigma)$. We consider several categories of models of $T$:

\begin{enumerate}
    \item $\Modhom(T)$ is the category of models of $T$ with homomorphisms for morphisms.
    \item $\Modstr(T)$ is the category of models of $T$ with strong embeddings for morphisms.
    \item For each pair of regular cardinals $\lambda \leq \kappa$, $\Modelem {\kappa,\lambda} (T)$ is the category of models of $T$ with $L_{\kappa,\lambda}$-elementary embeddings for morphisms.
\end{enumerate}
If $T$ is the empty theory, we write $\Str(\Sigma)$ for $\Modstr(\Sigma)$.
\end{Def}

The significance of basic theories for us is as follows.

\begin{fact}
By \cite[Theorem 5.42]{lpacc}, $\Modelem {\kappa,\lambda} \lambda(T)$ is accessible if $T$ is a theory of $L_{\kappa,\lambda}$. The categories $\Modhom(T)$ and $\Modstr(T)$ are \emph{not} always accessible, even when $T$ is in $L_\omega$. However, by \cite[Theorem 5.35]{lpacc}, $\Modhom(T)$ is accessible when $T$ is a basic theory. Moreover, if $T$ is a theory, then there is a theory $T'$ in an extended language with the same models as $T$ such that a homomorphism of $T'$-models corresponds to a strong embedding of $T$-models. Namely, we add a relation $R'$ for each relation symbol $R$ in the language of $T$ and an axiom $\forall x (R'(x) \leftrightarrow \neg R(x))$ and a relation symbol $R'_=$ such that $\forall x\forall y(R'_=(x,y) \leftrightarrow x\neq y)$. Note that the formula $\forall x (R'(x) \leftrightarrow \neg R(x))$ can be broken into $\forall x (R'(x)\rightarrow \neg R(x))$ and $\forall x (\neg R(x)\rightarrow R'(x))$, which is equivalent to $(\exists x(R'(x)\wedge R(x)))\rightarrow \false$ and $\forall x(\true\rightarrow R(x)\vee R'(x))$ respectively. The first formula is basic since it is an implication between positive-primitive formulae. The second one is again basic since $\true$ and $R(x)\vee R'(x)$ are positive-existential.
Hence if $T$ is basic, $T'$ is logically equivalent to a basic theory. Thus $\Modstr(T)$ is also accessible when $T$ is basic. 

In particular, since the empty theory is universal-positive, it is basic, and so $\Str(\Sigma)$ is accessible for any signature $\Sigma$.
\end{fact}

\begin{Def}[{\cite[Theorem 5.5]{beke-rosicky}}]\label{def:aec}
For a given functor $F:\calD \to \calC$, we say $F$  is \defterm{coherent} if for any objects $A,B,C \in \calD$ and $f:B\to C$ and $h:A\to C$. If there is $g' \in \Hom_\calC(FA,FB)$ such that $Fh=Ffg'$, then there is $g \in \Hom_\calD(A,B)$ such that $Fg=g'$. We say a subcategory $\calD$ of $\calC$ is \defterm{coherent} if the inclusion functor is coherent.\\
For a finitary signature $\Sigma$, an \defterm{Abstract Elementary Class (AEC)} in $\Sigma$ is a coherent, accessible subcategory of $\Str(\Sigma)$ that is closed under filtered colimits.\\
We say that $\calC$ has the \defterm{joint embedding property (JEP)} if for any finite, discrete diagram in $\calC$ there is a cocone. In other words, for every $A,B \in \calC$, there is an object $C$ and two morphisms $A\rightarrow C$, $B\rightarrow C$.\\
We say that $\calC$ has the \defterm{amalgamation property (AP)} if for any span in $\calC$, i.e. a diagram of the form $B \leftarrow A \to C$, there is a cocone. \\
We say that $\calC$ \defterm{has no maximal models} if, for every object $X \in \calC$, there is an object $Y$ in $\calC$ and a morphism $X \to Y$ which is not an isomorphism.
\end{Def}

\begin{eg}
A motivating example of an Abstract Elementary Class is the category $\Modelem {\kappa,\omega} (T)$, especially when $\kappa = \omega$. It is easy to see that $\Modelem\omega(T)$ has JEP and AP, and no maximal models. In \cite{classify}, we have shown that the classifying space of $\Modelem\omega(T)$ has fundamental group $\Gal_L(T)$, the so called Lascar group of $T$. It is unclear if $|\Modelem\omega(T)|$ can have non-vanishing higher homotopy groups.
\end{eg}

\begin{eg}[Axiomatizing Presheaves]\label{eg:ax-psh}
For $\calC$ a small category, there is a canonical finitary signature $\Sigma_\Psh(\calC)$ which has one sort $[C]$ for each object $C \in \calC$ and one unary function symbol $[f]: [D] \to [C]$ for each morphism $f: C \to D$ in $\calC$. There are no relation symbols. A $\Sigma_\Psh(C)$-structure $X$ comprises a set $X_C$ for each $C \in \calC$ and a function $X_f: X_D \to X_C$ for each $f: C \to D$ in $\calC$. There is a canonical theory $T_\Psh(C)$ in the language of $\Sigma(C)$ comprising the axiom $\forall x ([\id_C](x) = x)$ for each $C \in \calC$ and variable of sort $[C]$, and the axiom $\forall x ([f]([g](x)) = [gf](x))$ for each pair of composable morphisms $C\xrightarrow f D \xrightarrow g E$ in $\calC$ and variable $x$ of sort $[E]$. Thus a model $X$ of $T_\Psh(\calC)$ satisfies the equations $X_{\id_C}(x) = x$ and $X_f(X_g(x)) = X_{gf}(x)$, i.e. $X$ is precisely a presheaf on $\calC$. In fact, the category $\Modhom(T_\Psh(\calC))$ is canonically isomorphic to the category $\Psh(\calC)$ of presheaves on $\calC$, and the category $\Modstr(T_\Psh(\calC))$ is canonically isomorphic to the category $\Psh(\calC)^\mono$ of presheaves on $\calC$ with injective natural transformations for morphisms. We will identify these categories with their images under these isomorphisms.

In particular, $\Modstr(T_\Psh(\calC))$ is an Abstract Elementary Class.
\end{eg}

\section{The construction}

Let $\calC$ be a small category of monomorphisms, and let $T_\Psh(\calC)$ be the theory whose models are presheaves over $\calC$, as in Example \ref{eg:ax-psh}. 

\begin{Def}[Axiomatizing functors into $\calC \downarrow X$]
Let $K$ be a finite category, and let $F: K \to \calC$ be a functor. Let $\vec x = (x_k)_{k \in \Ob K}$ be tuple of variables the language of $T_\Psh(\calC)$ where $x_k$ has sort $F(k)$. Let $\Lift_F(\vec x)$ be the formula $\wedge_{f \in \Mor K} x_{\dom f} = [f](x_{\cod f})$.

Let $i: K \to L$ be a functor between finite categories, and let $G: L \to \calC$ be a functor such that $F = Gi$. Let $\vec y = (y_l)_{l \in \Ob L}$ be a tuple of variables where $y_l$ has sort $G(l)$. Let $\Ext_{i,G}(\vec x, \vec y)$ be the formula $\wedge_{k \in \Ob K} \vec y_{i(k)} = \vec x_k$.
\end{Def}

For any presheaf $X$ on $\calC$, write $\pi: \calC \downarrow X\to \calC$ for the projection.

\begin{lem}\label{lem:meaning}
Let $K$ be a finite category, and let $F: K \to \calC$ be a functor. Let $\vec x = (x_k)_{k \in \Ob K}$ be tuple of elements of a $T_\Psh(\calC)$-structure $X$ where $x_k$ has sort $F(k)$. Then $\Lift_F(\vec x)$ is satisfied if and only if the assignment $\tilde F: k \mapsto x_k$ determines a functor $\tilde F: K \to \calC \downarrow X$ such that $\pi \tilde F = F$.

Let $i: K \to L$ be a functor between finite categories, and let $G: L \to \calC$ be a functor; set $F = Gi$. Suppose that $\Lift_F(\vec x)$ and $\Lift_G(\vec y)$ hold for tuples $\vec x, \vec y$ in a $T_\Psh(\calC)$-structure $X$, and let $\tilde F: K \to \calC \downarrow X, \tilde G: L \to \calC \downarrow X$ be the corresponding functors. Then $\Ext_{i,G}(\vec x, \vec y)$ is satisfied if and only if $\tilde F = \tilde G i$.
\end{lem}
\begin{proof}
This is clear.
\end{proof}

\begin{Def}
Fix a functor $F: \sd^k (\partial \Delta[n]) \to \calC$. For each $k\geq 1$ and $n \in \nats$, we define the following sentence in the language of $T_\Psh(\calC)$:
\[\Asph_{n,k, F} := \forall \vec x (\Lift_F(\vec x) \to \bigdoublevee_{l \in \nats} \bigdoublevee_{G(\sd^{k+l} i) = Fp^l}  \exists \vec y (\Lift_G(\vec y) \wedge \Ext_{\sd^{k+l} i,G}( (p^l)^\ast(\vec x), \vec y) ))\]
Here the inner $\bigdoublevee$ is over all functors $G: \sd^{k+l} \Delta[n] \to \calC$ such that $G(\sd^{k+l} i) = Fp^l$. The functor $i: \partial \Delta[n] \to \Delta[n]$ is the inclusion as above, and the functor $p^l: \sd^{k+l}(\partial \Delta[n]) \to \sd^k(\partial \Delta[n])$ is the collapse map as above. The tuple $(x_j)_{j \in \Ob \sd^k (\partial \Delta[n])}$ is a tuple of variables where $x_j$ has sort $F(j)$. The tuple $(p^l)^\ast(\vec x)$ is the tuple $(x_{p^l(j)})_{j \in \Ob \sd^{k+l} (\partial \Delta[n])}$. The tuple $(y_j)_{j \in \Ob \sd^{k+l} \Delta[n]}$ is a tuple of variables where $y_j$ has sort $G(j)$. 

We also define the following sentence in the language of $T_\Psh(\calC)$ for $\sigma: d \to c \in \Ob \calC$:
\[\Inj_\sigma := \forall x,y ( [\sigma](x) = [\sigma](y) \to x = y)\]
We define $T_\contr(\calC)$ to be the theory extending $T_\Psh(\calC)$ by each instance of $\Asph_{n, k, F}$ for $n \in \nats$, $k \geq 1$, and $F: \sd^k (\partial \Delta[n]) \to \calC$, and each instance of $\Inj_\sigma$ for $\sigma$ a morphism of $\calC$.
\end{Def}

\begin{rmk}
Note that the $\bigdoublevee$ term of the sentence $\Asph_{n,k,F}$ is empty if $F$ admits no extension $G$ after any number of subdivisions; in this case $\Asph_{n,k,F}$ says that there are no lifts $\tilde F$ of $F$. For example, if $\calC$ has nontrivial homotopy type, this happens whenever $F$ corresponds to a nonzero element of $\pi_{n-1}(\calC)$.
\end{rmk}

\begin{lem}\label{lem:constprops}
\begin{enumerate}
    \item $\Asph_{n,k,F}$ is a basic formula of $L_{\kappa^+, \omega}$, where $\kappa = \max(\aleph_0, \#\calC)$; $\Inj_\sigma$ is a basic formula of $L_{\omega, \omega}$.
    \item Each instance of $\Asph_{n,k,F}$ holds of a model $X$ if and only if $\calC \downarrow X$ is contractible, where the slice is taken in $\Psh(\calC)$.
    \item Each instance of $\Inj_\sigma$ holds of a model $X$ if and only if every morphism $c \to X$ in $\Psh(\calC)$ from a representable is a monomorphism.
    \item Each representable presheaf on $\calC$ is a model of $T_\contr(\calC)$.
\end{enumerate}
\end{lem}
\begin{proof}
$(1):$ Follows from definition.

$(2):$ This follows from Lemma \ref{lem:meaning} and Fact \ref{fact:contr}.

$(3):$ This is clear.

$(4):$ $\calC$ is a category of monomorphisms, so the $\Inj$ sentences hold by $(3)$. Moreover, each slice category $\calC \downarrow c$ has a terminal object and so is weakly contractible, so the $\Asph$ sentences hold by $(2)$.
\end{proof}

\begin{thm}\label{thm:mainthm}
The inclusion $\calC \to \Modstr(T_\contr(\calC))$ induces a homotopy equivalences of classifying spaces. Moreover, $\Modstr(T_\contr(\calC))$ is an Abstract Elementary Class with amalgamation. If no connected component of $\calC$ is a groupoid, then $\Modstr(T_\contr(\calC))$ has no maximal models.
\end{thm}
\begin{proof}
Note that $\Modstr(T_\contr(\calC))$ is the category of models of $T_\contr(\calC)$ with injective maps as morphisms. Let $X \in \Modstr(T_\contr(\calC))$. The formulae $\Asph_{n,k,F}$ ensure by Lemma \ref{lem:constprops} that $\calC \downarrow X$ has a contractible classifying space when the slice category is computed in the category of presheaves. The formulae $\Inj_\sigma$ ensure by Lemma \ref{lem:constprops} that it does not matter whether we compute the slice category in the category of presheaves, or in the category of presheaves with injective maps for morphisms. Thus Quillen's Theorem A \cite{quillen} applies, and the inclusion $\calC \to \Modstr(T_\contr(\calC))$ (which is a fully faithful functor by Lemma \ref{lem:constprops} and the fact that monomorphisms are preserved and reflected by the Yoneda embedding) induces a homotopy equivalence of classifying spaces.

Because $\Modstr(T_\contr(\calC))$ is axiomatized by basic sentences (Lemma \ref{lem:constprops}), it is an accessible category by \cite{lpacc}, Theorem 5.35. The construction $\nerve (\calC \downarrow (-))$ commutes with colimits and monomorphisms, and weakly contractible simplicial sets are closed under filtered colimits and pushouts along injections. Thus $\Modstr(T_\contr(\calC))$ is closed under filtered colimits and pushouts in the presheaf category. The former (along with accessiblity) implies that $\Modstr(T_\contr(\calC))$ is an AEC, since it is a full subcategory of an AEC, and the latter implies that $\Modstr(T_\contr(\calC))$ has amalgamation.

Now we assume that no connected component of $\calC$ is a groupoid. To see that $\Modstr(T_\contr(\calC))$ has no maximal models, assume for contradiction that $X$ is a maximal model. If $X$ is not also a minimal model, then there is a non-isomorphism $A \to X$, and then $X \to X \cup_A X$ is a non-isomorphism, a contradiction. Thus $X$ is also a minimal model, i.e. the connected component of $X$ in $\Modstr(T_\contr(\calC))$ is a groupoid. Because the inclusion $\calC \to \Modstr(T_\contr(\calC))$ is a weak homotopy equivalence, it is in particular surjective on connected components, and because the inclusion is fully faithful, there is a connected component of $\calC$ contained in the connected component of $X$. Thus a connected component of $\calC$ is a groupoid, contrary to hypothesis. 
\end{proof}

\begin{cor}
For any small homotopy type $X$, there is an AEC $\calA$ with amalgamation and no maximal models whose classifying space is of the homotopy type $X$. Moreover, if $X$ is connected, one may assume that $\calA$ has joint embedding property.
\end{cor}
\begin{proof}
By \cite{thomason}, every small homotopy type $X$ is realized by the classifying space of a small poset $\calC$. We may assume that no connected component of $\calC$ is a groupoid -- e.g. we may replace any such component with the poset $0 \to 1$ without changing the homotopy type. Then by Theorem \ref{thm:mainthm}, $\Modstr(T_\contr(\calC))$ is an AEC with the same homotopy type as $\calC$, which has amalgamation and no maximal models. It is not hard to see that if a category $\calD$ satisfies amalgamation property and has no maximal models, then each connected component of $\calD$ has joint embedding property. Hence the ``moreover" part of the statement follows.
\end{proof}

\begin{rmk}
Here is a variant on the construction of $\Modstr(T_\contr(\calC))$. Let $\tilde \calC \subseteq \Psh^\mono(\calC)$ be the smallest full subcategory with the following properties:
\begin{enumerate}
    \item Every representable presheaf is in $\tilde \calC$.
    \item For every span $X_1 \leftarrow X_0 \to X_2$ in $\tilde \calC$, the pushout $X_1 \cup_{X_0} X_2$ in $\Psh(\calC)$ is also in $\tilde \calC$.
    \item For every filtered diagram $(X_i)_{i \in I}$ in $\tilde \calC$, the colimit $\varinjlim_i X_i$ in $\Psh(\calC)$ is also in $\tilde \calC$.
\end{enumerate}
Note that $\Psh^\mono(\calC)$ is closed under filtered colimits in $\Psh(\calC)$, and so $\tilde \calC$ is closed under filtered colimits in both $\Psh^\mono(\calC)$ and $\Psh(\calC)$. However, if $X_1 \leftarrow X_0 \to X_2$ is a span in $\tilde \calC$, then although any two monomorphisms of presheaves $X_1 \to Y \leftarrow X_2$ agreeing on $X_0$ \emph{do} induce a morphism of presheaves $X_1 \cup_{X_0} X_2 \to Y$, this map need not be a monomorphism. Thus $\tilde \calC$ does not have pushouts. Note that $\tilde \calC$ is a full subcategory of $\Modstr(T_\contr(\calC))$. In examples, it seems to often be the case that $\tilde \calC = \Modstr(T_\contr(\calC))$. However, it seems unlikely that this holds in general.

Similar arguments as in Theorem \ref{thm:mainthm} imply that the Yoneda embedding $\calC \to \tilde \calC$ is a weak homotopy equivalence. Moreover $\tilde \calC$ has filtered colimits, amalgamation, no maximal models (if no connected component of $\calC$ is a groupoid), and the inclusion $\tilde \calC \to \Psh^\mono(\calC)$ is coherent. By \cite[Theorem 6.17]{lpacc}, $\tilde \calC$ is accessible if we assume that Vop\u enka's principle holds, so under Vop\u enka's principle, we have that $\tilde \calC$ is an alternate realization of all the properties of Theorem \ref{thm:mainthm}. We do not know whether $\tilde \calC$ is accessible without the assumption of Vop\u enka's principle.
\end{rmk}

\section{Examples}\label{sec:examples}
\begin{eg}[Directed trees]\label{eg:dirtree}
Let $\calC$ be the category $V^\to_\to E$. Then $\Psh(\calC)$ is the category of directed multigraphs and homomorphisms, $\Psh(\calC)^\mono$ is the category of directed multigraphs and injective homomorphisms. The representables $\calC \subset \Psh(\calC)^\mono$ consist of the graph $V$ consisting of a point with no edges the graph $E$ consisting of two vertices with a single edge from one to the other. Thus a directed multigraph $G$ satisfies the formulas $\Inj_\sigma$ if and only if it has no loops. The condition that $\calC \downarrow G$ be connected is equivalent to $G$ being a connected graph, and the condition that $\calC \downarrow G$ be simply-connected is equivalent to saying that $G$ is acyclic. Thus, every object of $\Modstr(T_\contr(\calC))$ is a directed tree, i.e. a simple, acyclic, undirected graphs plus an orientation for each edge. Conversely, for every such graph $G$, the slice category $\calC \downarrow G$ is weakly contractible. So $\Modstr(T_\contr(\calC))$ is the category of directed trees with injective orientation-preserving homomorphisms.
% In this case, we don't need to invoke Vop\u enka's Principle to see that $\Modstr(T_\contr(\calC))$ is accessible.

Note that $\pi_1(\calC) = \ints$, so likewise the fundamental group of the category $\Modstr(T_\contr(\calC))$ of directed trees is $\ints$. More precisely, $|\calC| \simeq S^1$, and so $|\Modstr(T_\contr(\calC))| \simeq S^1$ as well.
\end{eg}

In the next example, we use the fact that anodyne extensions (i.e. injective weak homotopy equivalences of simplicial sets) contain the horn inclusions and are closed under pushouts along arbitrary map, as well as colimits of chains.
\begin{eg}[Undirected trees]\label{eg:undirtree}
Let $\calC$ be the category of finite sets of cardinality 1 or 2 and injective maps. Then $\calC$ has two objects $V,E$ and three nonidentity morphisms: \begin{tikzcd} V \ar[r,shift left,"s"] \ar[r,shift right, "t" below] & E \ar[loop right,"\tau"]\end{tikzcd}. The full subcategory on $E$ is $\bbB C_2$. Note that in $\Cat$, we have $\calC = \bbB C_2 \cup_{\{E\}} \{V \xrightarrow s E\}$. This motivates considering the simplicial set $B C_2^+ := BC_2 \cup_{\{E\}} \{V \xrightarrow s E\}$ where $BC_2$ is the nerve of $\bbB C_2$ and the pushout is now taken in $\sSet$. Note that the inclusion $B C_2 \to B C_2^+$ is an anodyne extension. An induction on the cell structure reveals that the inclusion $B C_2^+ \to \nerve \calC$ is also anodyne. Thus $|\calC| \simeq B C_2$.

The category $\Psh(\calC)$ consists of directed multigraphs equipped with an orientation-reversing involution on edges, and homomorphisms preserving the involution. The category $\Psh(\calC)^\mono$ is the same with just injective homomorphisms. The full subcategory of representables $\calC$ consists of the one-vertex graph $V$ with no loops and the graph $E$ with two vertices $v_1,v_2$ and two edges $v_1 \to v_2$ and $v_2 \to v_1$, exchanged by the involution $\tau$. A graph $G \in \Psh(\calC)$ satisfies the formulas $\Inj_\sigma$ if and only if there are no loops, and thus is equivalent to the data of a loop-free, undirected multigraph. Arguing as in Example \ref{eg:dirtree}, the category $\Modstr(T_\contr(\calC))$ consists of all trees: simple, undirected, connected, acyclic graphs.
% Again we don't need Vop\u enka's Principle to see that $\Modstr(T_\contr(\calC))$ is accessible.
Thus the classifying space of the category of trees is $B C_2$.
\end{eg}

\begin{eg}[Suspensions]\label{eg:susp}
Let $\calC_0$ be a small category of monomorphisms, let $\calC_0^\triangleleft$ be $\calC$ with an initial object freely adjoined, and let $\calC = \Sigma \calC_0 := \calC_0^\triangleleft \cup_{\calC_0} \calC_0^\triangleleft$. Then $\calC$ is also a small category of monomorphisms,\footnote{Note it's important here that we adjoined \emph{initial} objects. If we adjoined \emph{terminal} objects, or perhaps an initial and a terminal object, we would generally destroy the monomorphism property.} with $|\calC| \simeq \Sigma |\calC_0|$, where now $\Sigma X$ denotes the suspension of a space $X$. Any object $X$ of $\Psh(\calC)$ has a reduct $X_0$ which is an object of $\Psh(\calC_0)$; we call this the \defterm{underlying $\calC_0$-structure of $X$}. In addition, for any $X \in \Psh(\calC)$, there are two sets $X_1 = X(i_1)$, $X_2 = X(i_2)$ where $i_1,i_2$ are initial objects of the two copies of $\calC_0^\triangleleft$. We refer to $X_1$ (resp. $X_2$) as the \emph{first (resp. second) palette} of $X$. The data of $X \in \Psh(\calC)$ endows $X_0$ with certain maps to $X_1$ (resp. $X_2$) which we refer to as the \defterm{first (resp.) second coloring} of $X$. More precisely, if we view $X_0$ as a model of $T_\Psh(\calC_0)$, then the first and second coloring assign to each $x \in X_0$ a first color chosen from the $X_1$ palette and a second color chosen from the $X_2$ palette. The functoriality of $X$ as a presheaf ensures that the first (resp. second) color of $x$ depends only on its \defterm{connected component} in $X_0$. Here the connected components of $X_0$ are the components of the (essentially unique) maximal decomposition of $X_0$ as a coproduct of presheaves on $\calC_0$. In fact, the data of $X \in \Psh(\calC)$ is entirely encoded in the data of $X_0,X_1,X_2$ and the two colorings of $X_0$. Thus we may equivalently think of a presheaf on $\calC$ as a presheaf on $\calC_0$ equipped with an ordered pair of colorings of its connected components. Morphisms of $\Psh(\calC)$ are morphisms of $\Psh(\calC_0)$ equipped with maps of palettes which respect the colorings. In $\Psh(\calC)^\mono$, the morphisms of $\Psh(\calC_0)$ and the maps of palettes are restricted to be injective. An object of $\Psh(\calC)$ satisfies the sentences $\Inj_\sigma$ if and only if its underlying $\Psh(\calC_0)$-structure does. By Theorem \ref{thm:mainthm}, $\Modstr(\calC)$ has the homotopy type of $\Sigma |\calC|$.
\end{eg}

Unfortunately, in Example \ref{eg:susp} the contractibility condition on objects of $\Modstr(\calC)$ is difficult to understand in terms of $\calC_0$-structures. We will now construct a modified theory whose models are easier to understand in terms of $\calC_0$-structures.

\begin{Def}
Let $\calC_0$ be a small category of monomorphisms, and let $\calC = \calC_0^\triangleleft \cup_{\calC_0} \calC_0^\triangleleft$ as in Example \ref{eg:susp}. Let $T_\contr(\calC,\calC_0)$ comprise the axioms of $T_\contr(\calC)$, along with the axioms $\Asph_{n,k,F}$ for $n \geq 1$ applied to the reduct $X_0$ of $X$ to a presheaf on $\calC_0$.
\end{Def}

\begin{lem}
Let $\calC_0$ be a small category of monomorphisms, and let $\calC = \calC_0^\triangleleft \cup_{\calC_0} \calC_0^\triangleleft$ as in Example \ref{eg:susp}. A presheaf $X \in \Psh(\calC)$ is a model of $T_\contr(\calC,\calC_0)$ if and only if the following hold:
\begin{enumerate}
    \item $X$ is a model $T_\contr(\calC)$;
    \item The reduct $X_0$ of $X$ to a presheaf on $\calC_0$ is a disjoint union of models of $T_\contr(\calC_0)$.
\end{enumerate}
\end{lem}
\begin{proof}
Note first that the sentences $\Inj_\sigma$ for $\calC$ include the sentences $\Inj_\sigma$ for $\calC_0$. Then note that the axioms $\Asph_{n,k,F}$ for the reduct $X_0$, for $n \geq 1$, are equivalent to the condition that $\calC_0 \downarrow X$ has vanishing homotopy groups $\pi_n$ for $n \geq 1$. In other words, they say that each connected component of $\calC_0 \downarrow X$ is weakly contractible. Equivalently, the connected components $X_0^i$ of $X_0$ satisfy the condition that $\calC_0 \downarrow X_0^i$ is weakly contractible. Equivalently, the connected components $X_0^i$ of $X_0$ are models of $T_\contr(\calC_0)$.
\end{proof}

\begin{Def}
Let $\calC_0$ be a small category of monomorphisms, and let $\calC = \calC_0^\triangleleft \cup_{\calC_0} \calC_0^\triangleleft$ as in Example \ref{eg:susp}, with reduct $X_0$ to a $\calC_0$-presheaf and palettes $X_1,X_2$. The \defterm{coloring graph} $G(X)$ of $X$ is defined to be the following undirected bipartite multigraph on $X_1 \amalg X_2$. There is an edge from $x \in X_1$ to $y \in X_2$ for each connected component $X_0^i$ of $X_0$ whose first color is $x$ and second color is $y$.
\end{Def}

\begin{lem}
Let $X \in \Psh(\calC)$ be such that each connected component of the reduct $X_0$ to a presheaf on $\calC_0$ is a disjoint union of models of $T_\contr(\calC_0)$. Then $X$ is a model of $T_\contr(\calC,\calC_0)$ if and only if the coloring graph $G(X)$ is a tree.
\end{lem}
\begin{proof}
We have a pushout diagram of simplicial sets
\begin{equation*}
    \begin{tikzcd}
        \nerve(\calC_0 \downarrow X) \ar[r] \ar[d] & \nerve((\calC_0^\triangleleft)_1 \downarrow X) \ar[d] \\
        \nerve((\calC_0^\triangleleft)_2 \downarrow X) \ar[r] & \nerve(\calC \downarrow X)
    \end{tikzcd}
\end{equation*}
in which each morphism is a monomorphism. Thus this is a homotopy pushout. The category $(\calC_0^\triangleleft)_i \downarrow X$ deformation retracts onto the $i$th palette $X_i$, and thus is homotopically discrete. By hypothesis, $\nerve(\calC_0 \downarrow X)$ is also homotopically discrete. Thus the homotopy pushout $\nerve(\calC \downarrow X)$ is homotopy equivalent to the graph $G(X)$, and so is contractible if and only if $G(X)$ is a tree.
\end{proof}

We note that the data of a bipartite graph which is a tree, with the two ``parts" of the bipartition labeled $0,1$, is equivalent to the data of a tree with its vertices each labeled either $0$ or $1$, subject to the rule that adjacent vertices have different labels. Observe that there are exactly two ways to so label any tree. We call such an object a \defterm{$0,1$-vertex-labeled tree}.

\begin{cor}
Let $\calC_0$ be a small category of monomorphisms, and let $\calC = \calC_0^\triangleleft \cup_{\calC_0} \calC_0^\triangleleft$ as in Example \ref{eg:susp}. Then $\Modstr(\calC)$ is equivalent to the following category:
\begin{itemize}
    \item An object consists of a $0,1$-vertex-labeled tree equipped with a labeling of each edge by an object of $\Modstr(\calC_0)$;
    \item A morphism consists of an injective map of trees preserving the $0,1$-labeling, along with, for each edge, a morphism of $\Modstr(\calC_0)$ between the corresponding structures.
\end{itemize}
Moreover, the homotopy type of $\Modstr(\calC)$ is the suspension of the homotopy type of $\Modstr(\calC_0)$.
\end{cor}

\begin{eg}[An explicit AEC with homotopy type $S^2$]\label{eg:s2}
Let $\calC_0$ be the category $V^\to_\to E$ as in Example \ref{eg:dirtree}, and let $\calC = \calC_0^\triangleleft \cup_{\calC_0} \calC_0^\triangleleft$. Then $\Modstr(\calC,\calC_0)$ has the homotopy type of $S^2$. An object of this category is a $0,1$-vertex-labeled tree with each edge labeled by a directed tree. A morphism is an injection of $0,1$-vertex-labeled trees along with, for each edge, an injection of trees between the labels.
\end{eg}
\paragraph{Acknowledgements} The authors would like to thank Will Boney and Greg Cousins for helpful comments. TC was supported in part by NSF grant DMS-1547292.
\bibliographystyle{plain}
\bibliography{bibliography.bib}

\begin{thebibliography}{1}

\bibitem{lpacc}
Ji{\v{r}}{\'i} Ad{\'a}mek and Ji{\v{r}}{\'i} Rosick{\'y}.
\newblock {\em Locally presentable and accessible categories}, volume 189.
\newblock Cambridge University Press, 1994.

\bibitem{beke-rosicky}
Tibor Beke and Jir{\'i} Rosick{\'y}.
\newblock Abstract elementary classes and accessible categories.
\newblock {\em Annals of Pure and Applied Logic}, 163(12):2008--2017, 2012.

\bibitem{bourke}
John Bourke.
\newblock Equipping weak equivalences with algebraic structure.
\newblock {\em Mathematische Zeitschrift}, pages 1--25, 2017.

\bibitem{classify}
Tim Campion, Greg Cousins, and Jinhe Ye.
\newblock Classifying spaces and the lascar group.
\newblock {\em arXiv preprint arXiv:1808.04915}, 2018.

\bibitem{goerss_jardine}
Paul Goerss and J.F. Jardine.
\newblock {\em Simplicial homotopy theory}.
\newblock Springer Science \& Business Media, 2009.

\bibitem{quillen}
Daniel Quillen.
\newblock Higher algebraic {K}-theory: {I}.
\newblock In {\em Higher {K}-theories}, pages 85--147. Springer, 1973.

\bibitem{thomason}
R.~W Thomason.
\newblock Cat as a closed model category.
\newblock {\em Cahiers Topologie G{\'e}om. Diff{\'e}rentielle}, 21(3):305--324,
  1980.

\end{thebibliography}

\end{document}